\documentclass{amsart}
\usepackage{amssymb}
\usepackage{latexsym}
\usepackage{longtable}
\usepackage[greek,american]{babel}   
\usepackage[iso-8859-7]{inputenc}      %

\newtheorem{teor}{Theorem}[section]
\newtheorem{cor}[teor]{Corollary}
\newtheorem{lema}[teor]{Lemma}
\newtheorem{prop}[teor]{Proposition}

\newtheorem{defn}[teor]{Definition}
\newtheorem{obswr}[teor]{Observation}
\newtheorem{remarkwr}[teor]{Remark}

\newtheorem{examplewr}[teor]{Example}

\newenvironment{remark}{\begin{remarkwr}\begin{upshape}}{\end{upshape}\end{remarkwr}}

\newcommand{\Q}{\mathbb{Q}}
\newcommand{\Z}{\mathbb{Z}}
\newcommand{\F}{\mathbb{F}}

\newcommand{\C}{\mathbb{C}}

\newcommand{\Ql}{\Q_{\ell }}
\newcommand{\GQ}{G_{\Q}}

\newcommand{\Gal}{\operatorname{Gal\,}}
\newcommand{\GL}{\operatorname{GL}}

\newcommand{\End}{\operatorname{End}}

\newcommand{\disc}{\operatorname{disc}}

\newcommand{\Aut}{\operatorname{Aut}}
\newcommand{\Fr}{\operatorname{Fr}}

\newfont{\gotip}{eufb10 at 12pt}

\newcommand{\cO}{{\mathcal O}}

\newcommand{\om}{{\omega }}
\newcommand{\ra}{{\rightarrow }}
\newcommand{\lra}{{\longrightarrow }}
\newcommand{\qbar }{{\bar \Q }}
\newcommand{\Norm}{\mathrm{N}}

\newcommand{\n}{{\mathrm{n}}}

\newcommand{\R}{{\mathbb R}}
\newcommand{\M}{{\mathrm{M}}}

\include{thebibliography}

\begin{document}

\title[Which quaternion algebras act on a modular abelian variety?]{Which quaternion algebras act on a modular abelian variety?}

\author{Victor Rotger}\footnote{Partially supported by Ministerio
de Ciencia y Tecnolog\'{\i}a BFM2003-06768-C02-02}

\address{Universitat Polit\`{e}cnica de Catalunya,
Departament de Matem\`{a}tica Aplicada IV (EPSEVG), Av.\ Victor
Balaguer s/n, 08800 Vilanova i la Geltr\'{u} (Barcelona), Spain.}

\email{vrotger@ma4.upc.edu}

\subjclass{11G18, 14G35}

\keywords{Abelian variety, modular, endomorphism, Galois representation}

\maketitle

\textgreek{Gia es'ena mo'urh, pou se gn'wrhsa xafnik'a 'ena fjin'oporo sthn 'Anw P'olh ths Jessalon'ikhs.}

\begin{abstract} Let $A/\Q $ be a modular abelian
variety. We establish criteria to prevent a given quaternion algebra
$B$ over a totally real number field to be the endomorphism algebra
of $A$ over $\qbar $. We accomplish this by analyzing the
representation of $\Gal(\qbar /\Q )$ on the points of $N$-torsion of
$A$ at primes $N$ which ramify in $B$.
\end{abstract}

\section{Introduction}\label{intro}

Let $A/\Q $ be an abelian variety of $\GL_2$-type over $\Q $ of
dimension $g\ge 1$, by which we mean that the ring of endomorphisms
$R=\End_{\Q }(A)$ is an order in a number field $E$ of degree $[E:\Q
]=g$.

Let $K/\Q $ be the minimal number field over which {\em all}
endomorphisms of $A\times \qbar $ are defined and let $\cO
=\End_{K}(A)$. Let $B=\End_K^0(A):=\End_K(A)\otimes _{\Z } \Q $.

As a consequence of the recent work on Serre's modularity conjecture
by Khare-Wintenberger \cite{KaWi}, Dieulefait \cite{Di} and Kisin
\cite{Ki}, the generalized Shimura-Taniyama Conjecture holds true.
This amounts to say that $A$ is isogenous over $\Q $ to a factor of
the Jacobian variety $J_1(L)$ of the modular curve $X_1(L)$ for some
positive integer $L\geq 1$.

As in \cite[\S 1]{BFGR}, christen $(\cO , R)$ a {\em modular pair}.
The {\em minimal level} of $(\cO, R)$ is the minimal value of $L\geq
1$ as above. Similarly, we shall also say that a triplet $(\cO , R,
K)$ is modular if $\End_{\Q }(A)\simeq R$ and $\End_{K}(A)=\cO $ for
some modular abelian variety $A$.

A conjecture attributed to Coleman (cf.\,\cite[\S 1]{BFGR}) predicts
that, for any fixed $g\geq 1$, there exist only finitely many
isomorphism classes of modular pairs $(\cO , R)$. Also, the
possibilities for modular triplets $(\cO , R, K)$ should also be
very limited. The reader may consult \cite[\S 1]{BFGR} for further
motivation of this conjecture.

It is the purpose of these notes to address this question for {\em
absolutely simple modular abelian varieties $A/\Q $ such that $B
\supsetneq E$}.

Under this assumption, it follows from Albert's classification of
involuting simple algebras and the work of Shimura
that $g$ is even and $B$ is a
totally indefinite division quaternion algebra over a totally real
field $F$ of degree $[F:\Q ]=g/2$. In
particular, $E$ is a quadratic extension over $F$ and a maximal
subfield of $B$.

It is also known that $E$ is totally real and $K$ is imaginary
quadratic (cf.\,Lemma \ref{ole} (i), (ii)). Write $E=F(\sqrt{m})$
for some totally positive square-free integral element $m\in R_F$.
Although we choose $m$ square-free, note that there may exist a non
principal ideal $\mathfrak m_0$ of $F$ such that $\mathfrak
m_0^2\mid m R_F$.

Write also $K=\Q(\sqrt{-d})$ for some square-free integer $d\geq 1$.
By a theorem of Ribet (cf.\,\cite[Theorem 5.6]{Ri}), the values of
$m$ and $d$ are related by the isomorphism
\begin{equation}\label{Ribet}
B\simeq (\frac{-d, m}F).
\end{equation}

Write $\mathfrak D = \disc (B)$ for the reduced discriminant of $B$,
that is, the {\em square-free} product of the finite ideals of $F$
which ramify in $B$. Since $(\frac{-d, m}F)$ is split at all prime
ideals $\mathfrak N \nmid 2 d m$, it follows from (\ref{Ribet}) that
\begin{equation}\label{ad}
\mathfrak N \mid d \text{ for all } \mathfrak N \mid \mathfrak D,
\mathfrak N \nmid 2 m.
\end{equation}

In order to have a numerical flavour of the issue, let us report on
some explicit computations. The table below lists relevant numerical
data for all modular pairs $(\cO , R)$ of minimal level $L\leq 5400$
such that $\mathfrak D \ne (1)$ and $[F:\Q ] \leq 4$.

We refer the reader to \cite[Prop.\,1.3]{BFGR} for $L\leq 7000$ when
$F=\Q $; see also \cite{Ha} for explicit details when $L\leq 3000$
and $F=\Q $. When $[F:\Q ] = 2, 3$ or $4$, J.\ Quer performed the
necessary computations.

Column $L$ lists the minimal level of the modular pair. In column
$\mathfrak D$ we provide the norms $\mathrm{N}_{F/\Q }(\mathfrak N)$
of the prime ideals $\mathfrak N\mid \mathfrak D$.

$$
\begin{array}{|c|c|c|c|c|c|}
\hline L & [F:\Q ] & \disc(F) & \mathfrak D & \mathrm{N}_{F/\Q }(m) & \disc(K) \\
\hline  675   & 1 & 1 & [2, 3] & 2 & -3 \\
\hline  1568 & 1 & 1 & [2, 3] & 3 & -4 \\
\hline 243 & 1 & 1 & [2, 3] & 6 &  -3 \\
\hline 2700 & 1 & 1 & [2, 5] & 10 & -3 \\
\hline 1568 & 1 & 1 & [2, 7] & 7 &  -4 \\
\hline 3969 & 1 & 1 & [3, 5] & 15 & -7  \\
\hline
\hline  1089 & 2 & 5& [9, 11] & 11 & -3    \\
\hline 2592 & 2 & 33&  [2, 3] & 27 & -4 \\
\hline 3872 & 2& 5 & [4, 11] & 11 & -4 \\
\hline 3872& 2 & 5 & [4, 11] & 55 & -4 \\
\hline 4356& 2 & 5&  [5, 11] & 55& -3 \\
\hline 4761& 2 & 41&  [2, 5] & 10& -3\\
\hline
\hline 2187 & 3 & 81 & [3, 17] & 51 & -3 \\
\hline 2187 & 3 & 81 & [3, 8] & 24 &  -3 \\
\hline 3969 & 3 & 321 & [3, 3] & 81 & -7 \\
\hline 4563 & 3 & 1436 & [2, 3] & 6 & -3 \\
\hline
\hline 3267 & 4 & 5725 & [9, 11] & 11 & -3 \\
\hline 3267 & 4 & 13525 & [5, 9] & 5 & -3 \\
\hline
\end{array}
$$
\centerline{\ \ \ \ \ {\bf Table 1.} Modular pairs of minimal level $L\leq 5400$ and $[F:\Q ]\leq 4$.}
\medskip

\vskip 0.1 cm

\begin{defn}\label{Nl}
The set $\mathcal N_{\ell }$ of exceptional prime ideals of $F$ for
a given prime $\ell $ is
$$
\mathcal N_{\ell }\,=\, \{ \mathfrak N\,:\,\mathfrak N\mid \ell,\,\,
a^2 - s \ell \,\, \text { or }\,\,  a^4  - 4 a^2 \ell + \ell^2 \},
$$
for some $s= 0, 1, 2, 3, 4$ and some $a\in R_F$, $a\ne \sqrt{s
\ell}$, $\sqrt{2 \ell \pm \sqrt{3} \ell }$, such that $|\tau(a)|\leq
2 \sqrt{\ell } \text{ for all } \tau : F\hookrightarrow \R$.

\end{defn}

This is obviously a finite set of prime ideals, rather small for
small values of $[F:\Q ]$ and $\ell $. When $F=\Q $, we have for
instance $\mathcal N_{2} = \{ 2, 3, 5, 7\} $ and $\mathcal N_3 = \{
2, 3, 5, 11, 23\}$.

Given $\cO $, $R$, $K$, the results of this note provide necessary
conditions for the existence of a modular abelian variety $A/\Q $
such that $\End_{\Q }(A)\simeq R$ and $\End_K(A)\simeq \cO $. For
the sake of clarity and applicability, we gather in the theorem
below a simplified version of the results which are obtained
throughout the article. See the remaining sections for slightly more
general statements.

Below, let $\zeta_n$ denote a $n$-th primitive root of $1$. Also,
let $\cO_{\mathfrak D}$ denote a maximal order in the division
quaternion algebra $B$ of reduced discriminant $\mathfrak D\ne (1)$.
Recall that an order $\cO $ in $B$ is maximal if and only if $\disc
(\cO ) = \mathfrak D$ (cf.\,\cite[p.\,84]{Vi}).

\begin{teor}\label{main} Let $K$ be an imaginary quadratic field, let $F$ be a
totally real number field and let $m\in R_F$ be a square-free
totally positive element. If $(\cO _{\mathfrak D}, R_{F(\sqrt{m})},
K)$ is modular, then

\begin{enumerate}

\item [(i)] $m R_F = \mathfrak m_0^2\cdot \mathfrak m$, with $\mathfrak m_0$ and $\mathfrak m$ ideals of $R_F$, $\mathfrak m\mid \mathfrak D$. If $h(F)=1$, $\mathfrak m_0=1$. \vspace{0.1cm}

\item [(ii)] Any prime ideal $\mathfrak N\mid \mathfrak D$, $\mathfrak N\nmid 2m $, lies above a prime $N\equiv 3$ mod $4$. \vspace{0.1cm}

\item [(iii)] Assume $\mathfrak D\nmid 2 m$ and $\Q (\zeta_n+\zeta_n^{-1})\not\subset F$ for $n\ne 1, 2, 3, 4, 6$.
For any prime $\ell $ such that $(\frac{K}{\ell })\ne -1$ and
$\sqrt{\ell }$, $\sqrt{2 \ell }$, $\sqrt{3 \ell }$, $\sqrt{2\ell \pm
\sqrt{3} \ell }\,\not \in F$,  either

\begin{itemize}

\vspace{0.2cm}

\item $\mathfrak N\in \mathcal N_{\ell }$ for all $\mathfrak N\mid \frac{\mathfrak D}{(\mathfrak D, 2 m)}$, or

\vspace{0.2cm}

\item $(\frac{-\ell }{\mathfrak N }) \ne 1$ for all $\mathfrak N\mid \mathfrak
D$.




\end{itemize}

\end{enumerate}

\end{teor}

In the remaining sections of this note we shall develop the necessary machinery to prove this result.
Let us now illustrate and describe it in more detail.

Part $(i)$ is proved in Corollary \ref{m}. It was already known when $F=\Q $. In fact, if $F=\Q $, it was shown in \cite[Theorem 3.1]{RSY} that either
\begin{equation}\label{SY}
\mathfrak D = (m) \text{ or } (m N) \text{ for some prime } N.
\end{equation}

The proof of this result relies on a careful study of the closed
fibers of Morita's model of the Atkin-Lehner quotients of the
Shimura curve of discriminant $\mathfrak D$; it is reasonable to
expect that, under certain hypothesis, a similar result should hold
for higher degree number fields $F$. See the discussion following
Corollary \ref{kappa} for an approach to this question by a
different method.

In Table $1$ we quote the existence of a modular abelian variety
$A/\Q $ of dimension $4$ and level $3872$ such that $\End_{\Q }^0(A)
= F (\sqrt{ m})$, where $m\in F=\Q(\sqrt{5})$, $\mathrm{N}_{F/\Q
}(m) = 55$, and $\cO = \End_{\Q(\sqrt{-1})} (A)$ is an order in the
quaternion algebra $B$ over $F$ of discriminant $\mathfrak D =
2\cdot (4+\sqrt{5})$. By $(i)$, $\cO $ is not maximal in $B$.

Part $(ii)$ is proved in Corollary \ref{N} and $(iii)$ in Theorem
\ref{MaIn}.

As for $(iii)$, note firstly that for a given fixed $g$, the
condition on $F$ not to contain $\Q(\zeta_n+\zeta_n^{-1})$ for any
$n\ne 1, 2, 3, 4, 6$ excludes only finitely many fields $F$. For
instance, when $g=2$ this condition is empty, whereas for $g=4$ this
only excludes $F=\Q(\sqrt{2})$ and $\Q(\sqrt{5})$.

In fact, the restriction on $F$ can be relaxed (by allowing
$\zeta_n+\zeta_n^{-1}\in F$ for larger values of $n$) at the cost of
additional items in part $(iii)$. This can be achieved by following
the proof of Theorem \ref{MaIn}; we leave the details to the
interested reader.

The condition on $\ell $ not to be inert in $K$ can not be removed;
as for the requirement that $\sqrt{\ell }$, $\sqrt{2 \ell }$,
$\sqrt{3 \ell }$, $\sqrt{2\ell \pm \sqrt{3} \ell }\,\not \in F$, see
Theorem \ref{MaIn} for a statement without this hypothesis.

Roughly, $(iii)$ claims that for {\em half} of the primes $\ell $,
we have $(\frac{-\ell }{\mathfrak N}) = -1$ for all $\mathfrak N\mid
\mathfrak D$ unless all $\mathfrak N\mid \frac{\mathfrak D}{(\mathfrak D, 2
m)}$ lie in the exceptional set $\mathcal N_{\ell }$. Here, by half
of the primes we mean those $\ell $ such that $(\frac{K}{\ell })\ne
-1$ and $\sqrt{\ell }$, $\sqrt{2 \ell }$, $\sqrt{3 \ell }$, $\sqrt{2
\ell \pm \sqrt{3} \ell }\not \in F$.

When $F=\Q $ and $\mathfrak D = (m N)$ for an odd prime $N$, $(iii)$
-for $\ell =2, 3$, use the slightly stronger Theorem \ref{MaIn}-
says that for any $\ell $ such that $(\frac{K}{\ell }) \ne -1$,
either

\begin{itemize}
\item $N \in \mathcal N_{\ell }$, or
\item For all primes $p\mid \mathfrak D$: $(\frac{-\ell }{p})\ne 1$ if $\ell\ne 2$ (resp. $p\not \equiv 1$ mod $8$ if $\ell=2$).
\end{itemize}

When trying to apply $(iii)$ to preclude explicit families of pairs
$(\cO , R)$ from being modular, it naturally arises the question of
how does $K$ depend on $(\cO , R)$. By means of the theory of
descent applied to certain unramified covers of Shimura varieties,
we prove in \cite{Ro} the following fact. See \cite[\S 5]{RSY} for a
particular case with $F=\Q $.

\begin{prop}\label{unr} Let $(\cO _{\mathfrak D}, R_{F(\sqrt{m})}, K)$ be a modular triplet. Assume that

\begin{itemize}

\item $m R_F$ is a square-free ideal, with $m=3$ or $\tau (m)>4$ for some embedding $\tau : F \hookrightarrow \R $, and

\item There exists $\mathfrak N \mid \mathfrak D$, $\mathfrak N \nmid 2 m$, of {\em odd} residual degree $f(\mathfrak N/N)$.

\end{itemize}

Then $K$ is unramified away from $\Norm_{F/\Q }(\mathfrak D) \cdot \disc(F/\Q )$.
\end{prop}

The virtue of the above result is that, when the hypothesis is
satisfied, there is a {\em finite} number of possibilities for $K$.
Observe that, when $F=\Q $, the hypothesis of Proposition \ref{unr}
simply require that $\mathfrak D = (m N)$ with $m, N > 2$.

In addition, notice that the isomorphism (\ref{Ribet}) $B\simeq
(\frac{-d, m}F)$ and Corollary \ref{kappa} impose further
restrictions on $K$. Let us illustrate how the combination of these
results can be used to prevent triplets $(\cO , R, K)$ from being
modular.

\begin{teor}\label{Sh} Let $F=\Q $. Let $M, N$ be two
different odd primes and assume $(\cO_{M N}, R_{\Q(\sqrt{M})})$ is a
modular pair. Then $N\equiv 3$ mod $4$ and $(\frac{-N}M)=-1$.
Moreover,

\begin{enumerate}

\item[(i)] If $M\equiv 3$ mod $4$, then $K=\Q(\sqrt{-N})$ and
$(\frac{-\ell }M) = -1$ for any odd prime $\ell $ such that
$(\frac{\ell }{N }) = 1$ and $N\not \in \mathcal N_{\ell }$.

\item[(ii)] If $M\equiv 1$ mod $4$, then $K=\Q(\sqrt{-d})$ with
$d=N$ or $M N$. If $d=N$, then $(\frac{-\ell }M) = -1$ for any $\ell
$ such that $(\frac{\ell }{N }) = 1$ and $N\not \in \mathcal N_{\ell
}$. If $d=M N$, then $N\equiv 3$ mod $8$ and $N\in \mathcal N_{\ell
}$ for all odd primes $\ell $ such that $(\frac{-M N}{\ell }) = 1$.

\end{enumerate}

\end{teor}

This is proved in Section \ref{last}. Since the finite sets
$\mathcal N_{\ell }$ are easily computable, notice that, by taking
explicit values of $\ell $, the above criterion precludes infinitely
many pairs $(\cO_{M N}, R_{\Q(\sqrt{M})})$ from being modular.
Similar examples with higher degree number fields $F$ can easily be
worked out in the same way.

The material of this note can be used to prove the non-existence of rational points (or at least, non-trivial rational points) on Atkin-Lehner quotients of Shimura varieties associated with totally indefinite quaternion algebras over totally real number fields. Care must be taken though because rational points over a field $K$ on these varieties do not quite correspond to abelian varieties (with extra structure) defined over $K$, but rather to those with $K$ as field of moduli; details will appear in \cite{Ro}. 

See also \cite{PY} for (non-overlapping) results in this direction in dimension $1$, where it should be noticed that the methods used there in are completely different.

\vspace{0.2cm}

{\bf Acknowledgement.} It is a pleasure to thank B. Jordan, J. Gonzalez, D. Lorenzini, J. Quer and the referee for their helpful remarks.

\section{Quaternion multiplication over number fields}\label{canonical}

Let $A/\Q $ be an abelian variety of (even) dimension $g$ of
$\GL_2$-type over $\Q $ and let $R=\End_{\Q }(A) \subset E=\End_{\Q
}(A)\otimes \Q $. Let $K$ be the minimal field over which all
endomorphisms of $A\otimes \qbar $ are defined. Assume $\cO
=\End_{K}(A)$ is an order in a division totally indefinite
quaternion algebra $B$ over a totally real number field $F$.

For any prime $N$, the action of $\GQ$ on the Tate module $V_{N}(A)$ of $A$
induces a Galois representation
$$
\mathcal R_{N} : \GQ \, \lra\, \GL_{2 g}(\Q_{N}),
$$
which in fact takes values in $\GL_2(R_E)$, independently of $\mathfrak N$.

For any prime ideal $\mathfrak N$ of $E$ over $N$, let $E_{\mathfrak N}$ stand for the completion of $E$ along
$\mathfrak N$ and write $V_{\mathfrak N}(A) = V_N(A)\otimes_{\Q_N}E_{\mathfrak N}$. Since the representation $\mathcal R_N$
is compatible with the action of $E$, it induces a Galois representation
$$
r_{\mathfrak N}: \, \GQ \, \lra \, \GL_2(E_{\mathfrak N}).
$$

\begin{lema}\label{good}
Let $\ell $ be a prime number. There exists a
finite extension $L_{\ell }/\Ql$ such that the closed fibre of the N\'eron model of $A\otimes_{\Q } L_{\ell }$ over the
ring of integers of $L_{\ell }$ is an abelian variety $\tilde A$ over $\F_{\ell }$.
\end{lema}

{\em Proof.} By \cite[Theorem 3]{Ri81}, $A$ has potential good
reduction at $\ell $. Under this assumption, Serre and Tate explicitly
construct $L_{\ell }$ and $\tilde A/\F_{\ell }$ at the end of p.\,498 of \cite{SeTa}.
$\Box $

\begin{remark} Let $A/\Q_{\ell }$ be an abelian variety with potential good reduction. By definition, there exist finite extensions $L/\Q_{\ell }$ over which $A$ acquires good reduction. Lemma above claims that $L$ can be chosen to have residue field $\F_{\ell }$. Note that in general there does not exist a {\em minimal} extension $L/\Q_{\ell }$ over which $A$ has good reduction (cf.\,\cite[p.\,498]{SeTa} for a proof of this fact when the residue field is algebraically closed). Hence, there may exist good reduction extensions $L/\Q_{\ell }$ for $A$ with non trivial residual degree such that $A\times M$ has bad reduction for any subextension $M\varsubsetneq L$.
\end{remark}

Let $\ell $ be a prime and let $\tilde
A/\F_{\ell }$ denote the closed fibre of the N\'eron model of $A\otimes L_{\ell }$ as in the lemma above.

Let $\Fr_{\ell }\in \Gal(\bar \F_{\ell }/\F_{\ell })$ be the
Frobenius automorphism. Let $\varphi_{\ell }$ be an element of
$\Gal(\qbar_{\ell }/L_{\ell })$ whose image in $\Gal(\bar \F_{\ell
}/\F_{\ell })$ under the canonical reduction map is $\Fr_{\ell }$.
By Serre-Tate's criterion \cite[Thm. 1]{SeTa}, $r_{\mathfrak
N}(\varphi_{\ell })$ is a well-defined element of
$\GL_2(R_E)$ up to conjugation.

Let $P_{\ell }(T) :=
P_{\varphi_{\ell }}(T) \in R_E[T]$ denote the characteristic polynomial of $r_{\mathfrak N}(\varphi_{\ell })$. We shall write $a_{{\ell }} = \mathrm{Tr}(r_{\mathfrak N}(\varphi_{\ell }))\in R_E$. By the work of Weil, $|\tau (a_{\ell })|\leq 2 \sqrt{\ell }$ for any embedding $\tau : E \, \ra \, \C $. By \cite[Prop. 3.5]{Ri}, $E=\Q(\{a_{\ell }\}_{\ell \nmid N_A})$, where $N_A$ stands for the conductor of $A$.

There is a natural action of $G_{\Q }$ on the ring of endomorphisms $\cO =\End_{\qbar }(A)$ of $A\otimes \qbar $. By the Skolem-Noether theorem, for any $\sigma \in G_{\Q }$ the automorphism $B\,\ra \, B$, $\beta \,\mapsto \beta^{\sigma }$ is inner: there exists $\om_{\sigma } \in \cO $ such that $\beta^{\sigma }=\om_{\sigma } \beta \om_{\sigma }^{-1}$. Thus, $\om_{\sigma }$ belongs to the normalizer $\mathrm{Norm}_B(\cO )$ of $\cO $ in $B$, because $\beta^{\sigma }\in \cO $ for any $\beta \in \cO $.

Since the endomorphisms of $R\subset \cO $ are defined over $\Q $, it follows that $\om_{\sigma }$ belongs to the commutator of $E$ in $B$, which is $E$ itself because it is a maximal subfield of $B$. Hence, $\om_{\sigma }\in R=E\cap \cO $. This induces a continuous homomorphism
$$
\psi : \, G_{\Q } \, \longrightarrow \, E^*/F^*, \quad \sigma \,\mapsto \, \om_{\sigma }.
$$



\begin{lema}\label{ole}

\begin{enumerate}
\item[(i)] $K$ is imaginary quadratic, say $K=\Q (\sqrt{-d})$, $d>0$.
\item[(ii)] $E=F(\om )$ for an element $\om \in \cO \cap \mathrm{Norm}_B(\cO )$, $\om^2=m$, where $m\in R_F$ is square-free and totally positive.
\end{enumerate}
\end{lema}

{\em Proof. } By \cite[Lemma 3.1]{Ri}, $\det (r_{\mathfrak N}) = \epsilon \cdot \chi_{N} : \GQ \, \lra \, E^*$, where $\chi_N$ is the $N$-adic cyclotomic character and $\epsilon : G_{\Q }\,\ra \, E^*$ is a character of finite order.

Let us first show that $E$ is totally real. Indeed, if it were not,
$E$ would be a CM-field. Let $\mathfrak L$ be the set of primes
$\ell \nmid N_A$ such that $a_{\ell }\ne 0$. If $\epsilon
(\varphi_{\ell })=1$ then $a_{\ell }\in F^*$ by \cite[Proposition
3.4]{Ri} and the fact that $E$ is a totally imaginary quadratic
extension of $F$. Since $\psi (\varphi_{\ell }) \in a_{\ell }\cdot
F^*$ by \cite[Theorem 5.5]{Ri}, we deduce in turn that $\psi
(\varphi_{\ell })\in F^*$, thus $\varphi_{\ell }\in \mathrm{Ker
}(\psi )$. \v{C}ebotarev Density Theorem together with the fact that
$\mathfrak L$ has density $1$ within the set of all primes
(cf.\,\cite{SeIHES}) would imply that $\mathrm{Ker }(\epsilon
)\subseteq \mathrm{Ker }(\psi )$, whence $K\subseteq
\qbar^{\mathrm{Ker }(\epsilon )}$. By \cite[Lemma 3.2]{Ri},
$\epsilon(c)=1$ for any complex conjugation $c\in G_{\Q }$ on $\qbar
$ and we deduce that $K$ would be totally real. Theorem $2$ in
\cite{Ri81} would apply to say that $B$ is not division, which
contradicts our assumptions.

Since $E$ is totally real, it contains no roots of unity $\zeta \ne
\pm 1$ and any non-trivial finite subgroup of $E^*/F^*$ has order
$2$. As $\psi $ embeds $\mathrm{Gal}(K/\Q )$ in $E^*/F^*$, we obtain
that $K$ is quadratic. Theorem 2 of \cite{Ri81} implies that it is
imaginary and $(i)$ follows.

Write $\Gal(K/\Q ) = \{ 1, \sigma \}$ and let $\om \in R$ such that
$\psi (\sigma )\in \om \cdot F^*$. As we argued above, $\om \in \cO
\cap \mathrm{Norm}_B(\cO )$. Since $\sigma ^2=1$, $\om^2=m\in R_F$
and we can assume $m$ is square-free. This shows $(ii)$. $\Box $

\vspace{0.2cm}

Recall that an order $\cO $ in $B$ is an {\em Eichler order} if it
is the intersection of two maximal orders (cf.\,\cite[p.\,39,
84]{Vi}). All orders with {\em square-free} discriminant $\disc (\cO
)$ are Eichler orders by \cite[p.\,84 and p.\,98, Ex.\ 5.3]{Vi}.

\begin{cor}\label{m}
Assume $\cO $ is an Eichler order. Then $m\cdot R_F = \mathfrak
m_0^2\cdot \mathfrak m$, where $\mathfrak m_0$ and $\mathfrak m$ are
integral ideals of $F$ and $\mathfrak m\mid \disc (\cO )$. If the
class number of $F$ is $h(F)=1$, then $m\mid \disc (\cO )$.
\end{cor}

{\em Proof. } Let $\om \in \cO \cap \mathrm{Norm}_B(\cO )$ as in the
above proof. By \cite[p.\,99, Ex.\ 5.4]{Vi}, $\om $ normalizes $\cO
_{\wp }$ in $B_{\wp }$ for any prime ideal $\wp $ of $F$. By
\cite[p.\,34 and 40]{Vi}, there are no restrictions on $m=\n (\om )$
at any $\wp \mid \disc (\cO )$. By \cite[p.\,40]{Vi}, $m\cdot
R_{F_{\wp }}$ is an even power of $\wp R_{F_{\wp }}$ for any $\wp
\nmid \disc (\cO )$. This shows the first part. If $h(F)=1$, all
ideals of $R_F$ are principal. Since $m$ is square-free, $\mathfrak
m_0=(1)$.  $\Box $

\vspace{0.3cm}

Since $E$ is totally real, the character $\epsilon $ mentioned in
the proof above is trivial (cf.\,\cite[p.\,244]{Ri}) and
\cite[Theorem 5.3]{Ri} asserts that $F=\Q(\{a^2_{\ell }\}_{\ell
\nmid N_A})$. Also, by \cite[Lemma 3.1]{Ri} $\det (r_{\mathfrak
N}(\varphi_{\ell })) = \ell $  and $P_{\ell }(T) = T^2 -a_{\ell } T
+ \ell $ for any prime $\ell $.

For primes $\ell \nmid N_A$ such that $a_{\ell }\ne 0$ we have by
\cite[Theorem 5.5]{Ri} that if $\ell $ remains inert in $K$, then
$a_{\ell }\in F^{*}\cdot \om $. Otherwise, $a_{\ell }\in
F^{*}$ and $P_{\ell }\in R_F[T]$.

In the terminology of modular forms, this is equivalent to saying that $A$ is the abelian variety attached (up to isogeny) to a normalized newform $f = \sum a_n q^n \in S_2(\Gamma_0(L))$ with an {\em extra-twist} by $\psi $. As we have seen, $F=\Q (\{ a^2_{\ell }\})$, $E=F(\sqrt{m})=\Q (\{ a_{\ell }\})$, $K=\Q(\sqrt{-d})=\mathrm{Ker }\psi $ and $B=(\frac{-d, m}F)$ are completely determined by the coefficients $a_{\ell }$. With this perspective in mind, the conjecture addressed in this note claims that the Fourier coefficients of a normalized newform for $\Gamma_0(L)$, $L\geq 1$, generate only finitely many endomorphism algebras $E$, $B$ of given degree over $\Q $.

\section{The case $\mathfrak N$ inert in $E$}\label{Ass}

Keep the notation of Section 2. In addition,
for the rest of the article we shall assume the following

\vspace{0.1cm}

{\bf Assumption.} {\em There exists a prime ideal $\mathfrak N_0\mid \disc(B)$, $\mathfrak N_0\nmid 2 m$, at which $\cO $ is locally maximal.}

\vspace{0.1cm}

Fix such a prime $\mathfrak N_0$ and let $R_0=R\cap F$, a suborder
of the ring of integers $R_F$ of $F$. By saying that $\cO $ is
locally maximal at $\mathfrak N_0$ we mean that $R_{0, \mathfrak
N_0}$, the localization of $R_0$ at $\mathfrak N_0$, is the ring of
integers of $F_{\mathfrak N_0}$ and that $\cO_{\mathfrak N_0}$ is
the (single) maximal order of $B_{\mathfrak N_0}$. Let $k_0$ be the
residue field of $R_0$ at $\mathfrak N_0$, say of $\nu := N^f$
elements for $f=f_{\mathfrak N_0/N}\geq 1$.

Having this, the condition on $\mathfrak N_0$ not to ramify in $E$
actually implies that $\mathfrak N_0$ remains inert in $R$. This
follows from the fact that $R$ embeds optimally in $\cO $, as by
construction $R=E\cap \cO $ (cf.\,\cite[Ch.\,II, \S 3]{Vi} for
details on Eichler's theory of optimal embeddings). Let $k$ be the
residue field of $\mathfrak N=\mathfrak N_0\cdot R$, which is a
quadratic extension of $k_0$. Reducing mod $\mathfrak N$ we obtain
the residual representation
$$
\bar{r}_{\mathfrak N}:\,G_{\Q }\,\longrightarrow \, \GL_2(k),
$$
that is, the representation of $G_{\Q }$ on
$A[\mathfrak N]=\{ x\in A:\,\beta \cdot
x=0\quad \forall \beta \in \mathfrak N \, \}$ of $A$.

Let $\pi_0$ be an uniformizer of $F_{\mathfrak N_0}$. By
\cite[Ch.\,II, \S 1]{Vi}, locally at $\mathfrak N_0$ the division
quaternion algebra $B_{\mathfrak N_0}$ can be described as $B_{\mathfrak N_0} =
E_{\mathfrak N} + E_{\mathfrak N} \pi $, where $\pi^2=\pi_0$ and $\beta \cdot \pi =
\pi \cdot {\beta }^{\tau }$ for any $\beta \in E_{\mathfrak N}$. Here, $\tau $ denotes the
non-trivial automorphism of $\Gal(E_{\mathfrak N}/F_{\mathfrak N_0})$.
Moreover, $\cO_{\mathfrak N_0} = R_{\mathfrak N} + \cdot
R_{\mathfrak N} \cdot \pi $.

Let $I=\{ \beta \in \cO : \n(\beta )\in \mathfrak N_0\}$. This is an
ideal of $\cO $ which locally at $\mathfrak N_0$ is $I_{\mathfrak N_0} =
\mathfrak N\cdot R_{\mathfrak N} + R_{\mathfrak N}\cdot \pi
$, whereas $I_{\mathfrak M_0} = \cO_{\mathfrak M_0}$ at the remaining
primes $\mathfrak M_0\ne \mathfrak N_0$ of $F$.
Let $C=\mbox{ Ker }(I: A\,\ra \, A) = \bigcap _{\beta \in I} \mbox{
Ker } (\beta : A \,\ra \,A)$, a subgroup of $A[\mathfrak
N]$.

The action of $B$ on $A(\qbar )$ makes $V_{\mathfrak N}(A)$ a
$B_{\mathfrak N_0}$-module, which must be free because $B_{\mathfrak N_0}$ is simple.
In fact, since $\dim_{F_{\mathfrak N_0}}(B_{\mathfrak N_0}) = \dim (V_{\mathfrak N}(A)) = 4$, we have
$V_{\mathfrak N}(A)\simeq B_{\mathfrak N_0}$.

In the same fashion, $T_{\mathfrak N}(A)$ is naturally a
module over $\cO _{\mathfrak N_0}$. In fact, one can identify $T_{\mathfrak N}(A)$
as a left ideal of $\cO _{\mathfrak N_0}$ by choosing an isomorphism $V_{\mathfrak N}(A)\simeq B_{\mathfrak N_0}$.

Since $\cO_{\mathfrak N_0}$ is a maximal order, all its ideals are principal by \cite[p.\,34]{Vi} and
we may write $T_{\mathfrak N}(A)=\cO_{\mathfrak N_0}\cdot x$ for some $x\in
A[\mathfrak N]$.

Note that $A[\mathfrak N] = T_{\mathfrak N}(A)/\mathfrak N\cdot
T_{\mathfrak N}(A)\simeq \cO_{\mathfrak N_0}/\mathfrak N\cdot
\cO_{\mathfrak N_0}$. Also, we have $C = \cO_{\mathfrak
N_0}/I_{\mathfrak N}\simeq R_{\mathfrak N}/\mathfrak N\cdot
R_{\mathfrak N}\simeq k$. As a $R$-module, $\Aut_{R}(C)\simeq k^*$.

The action of $G_{\Q }$ on $\cO $
leaves $I$ invariant because it is the single two-sided ideal of
$\cO $ of reduced norm $\mathfrak N_0$. Since $R=\End_{\Q }(A)$,
$G_{\Q }$ acts $R$-linearly on $C$. It thus induces a Galois
representation
$$
\alpha_{\mathfrak N}\,:\,G_{\Q }\,\longrightarrow \, \Aut_{R}(C)=k^*.
$$

As in \cite{Jo}, the character $\alpha_{\mathfrak N}$ shall play a key role in what follows, since it
encodes many of the arithmetic properties of the Galois representation of $G_{\Q }$ on the torsion of $A$.
Recall that $\nu = \sharp k_0$ and that $\psi $ is a character of $G_{\Q }$ with values in $\{ \pm 1\}$.

\begin{lema}\label{Borel} There exists a $k$-basis of $A[{\mathfrak N}]$ with
respect to which
$$\begin{matrix}
\bar {r}_{\mathfrak N} : &  \GQ  & \lra & \GL_2(k) \\
    &   \sigma & \ra & \begin{pmatrix}
    \psi(\sigma ) \alpha_{\mathfrak N}(\sigma )^{\nu } & 0 \\ \beta_{\sigma } & \alpha_{\mathfrak N}(\sigma )
    \end{pmatrix}
\end{matrix}$$
for some $\beta_{\sigma }\in k$.
\end{lema}

{\em Proof. }  Write $\cO_{\mathfrak N_0}=R_{\mathfrak N}+R_{\mathfrak N}\cdot \pi $ as above and let
$x\in A[\mathfrak N]$ such that $A[{\mathfrak N}]=\cO_{\mathfrak N_0}/{\mathfrak N}\cO_{\mathfrak N_0}\cdot x = R_{\mathfrak N}/{\mathfrak N} R_{\mathfrak N}\cdot
x\,+\,R_{\mathfrak N}/{\mathfrak N} R_{\mathfrak N}\cdot \pi (x)$. We shall compute
$r_{\mathfrak N}(\sigma )$ for $\sigma \in G_{\Q }$ with respect to the $k$-basis $\{ x, \pi
(x)\}$ of $A[{\mathfrak N}]$.

Let $\sigma \in G_K$. Since all endomorphisms of $A$ are defined
over $K$ it follows that $\pi^{\sigma } = \pi $. Let $x^{\sigma } =
\alpha_{\sigma }\cdot x\,+\,\beta_{\sigma }\cdot \pi(x)$ for some
$\alpha_{\sigma }, \beta_{\sigma }\in R_{\mathfrak N}$. Then $(\pi(x))^{\sigma }
= \pi^{\sigma }(x^{\sigma }) = \pi (\alpha_{\sigma }\cdot
x\,+\,\beta_{\sigma }\cdot \pi(x)) = {\alpha}^{\tau }_{\sigma } \pi
(x) + {\beta}^{\tau}_{\sigma } \pi^2(x) = {\alpha}^{\tau}_{\sigma }
\pi (x)$, because $\pi^2=\pi_0$ and $\pi_0 \cdot x= 0$. Switching
$\alpha_{\sigma }$ by ${\alpha}^{\tau}_{\sigma }$ and reducing mod
${\mathfrak N}$ it follows that $r_{\mathfrak N}(\sigma ) = \begin{pmatrix} \alpha_{\sigma
}^{\nu } & 0 \\ \beta_{\sigma } & \alpha_{\sigma }\end{pmatrix}$. In
particular we recover the canonical character $\alpha_{\mathfrak N} : G_{\Q
}\,\ra \, k^*$ as $\alpha_{\mathfrak N}(\sigma )=\alpha_{\sigma }$.

Let $c \in G_{\Q }$ denote a complex conjugation. As in Lemma \ref{ole} and its proof, for all $\beta \in \cO $ we have $\beta^c=\om
\beta \om^{-1}$, $\om \in R$, $\om^2=m\in F^*$.

It follows that $\pi^c = \om \pi \om^{-1} = \pi {\om }^{\tau}
\om^{-1} = -\pi $.  Write $x^c = \alpha_c x + \beta_c \pi (x)$ for
some $\alpha_c, \beta_c \in R_{\mathfrak N}$. Then $(\pi (x))^c = \pi^c(x^c) =
-\pi (\alpha_c x + \beta_c \pi (x)) = -\pi_0 {\beta}^{\tau}_c x -
{\alpha}^{\tau}_c \pi (x)  = -{\alpha}^{\tau}_c \pi (x)$. In
particular $\alpha_{\mathfrak N}(c)=-\alpha^{\tau}_c$. Since $c^2=\mbox{Id}$,
$\alpha^{\tau }_c=\pm 1$ and $\bar r_{\mathfrak N}(c) =
\begin{pmatrix} \pm 1& 0 \\ \beta_c & \mp 1\end{pmatrix}$. Either value takes $\bar r_{\mathfrak N}(c)$, the lemma
follows. $\Box $

\vspace{0.1cm}

Let $\GQ^{ab} = \Gal(\Q^{ab}/\Q )$ denote the Galois group of the
maximal abelian extension $\Q^{ab }$ of $\Q $. For every prime $\ell
$, let $\Q_{\ell}^{ab}$ denote the maximal abelian extension of
$\Q_{\ell }$ and let $G_{\ell }^{ab }=\Gal(\Q_{\ell}^{ab}/\Q_{\ell
})$. Let $I_{\ell }^{ab}$  denote the inertia subgroup of $G_{\ell
}^{ab}$, which we regard as a subgroup of $G_{\Q }^{ab}$. Local
class field theory establishes a canonical isomorphism
$$
\omega_{\ell }:\, \Z_{\ell }^* \,\stackrel{\sim}{\lra }\, I_{\ell }^{ab},
$$
the local Artin reciprocity map. Global class field theory yields a surjective map
$$
\omega : \, \prod_{\ell }\Z_{\ell }^* \, \stackrel{\prod \om_{\ell }}{\twoheadrightarrow } \, \GQ^{ab}.
$$

We shall denote $$\varrho _{\ell } : \,\Z_{\ell }^* \,
\stackrel{\om_{\ell }}{\simeq }\, I_{\ell }^{ab }\,\subset \,
\GQ^{ab }\, \stackrel{\alpha_{\mathfrak N}}{\lra }k^*$$ the
composition of the local Artin reciprocity map together with the
character $\alpha_{\mathfrak N}$ acting on the canonical torsion
subgroup $C$. In particular, since $\varrho _N: \,\Z_N^* \,\lra \,
k^*$ is a continuous homomorphism, $\varrho _N (\Z_N^*)\subseteq
\F_N^*\,\subset \, k^*$.

\begin{defn}
Let $\kappa (A, \mathfrak N)$ be the least common multiple of
$4$ and the orders of $\alpha_{\mathfrak N}(I_{\ell })$ as $\ell $
runs over all prime integers $\ell \ne N$.
\end{defn}

Since $A$ has good reduction outside a finite set of primes, and
potential good reduction at those, $\kappa (A, \mathfrak N)$ is a well-defined
integer.

For a prime $\ell $, let $\mathfrak a_{\ell
}= \prod \mathfrak L$, where $\mathfrak L$ runs over the prime
ideals of $F$ above $\ell $ such that $[F_{\mathfrak L}: \Q_{\ell
}]$ is odd. Let $B(\ell )$ be the totally definite quaternion algebra
over $F$ of discriminant $\frac{\mathfrak D\cdot \mathfrak a_{\ell }}
{(\mathfrak D, \mathfrak a_{\ell })^2}$.

\begin{defn}\label{B}
Let $\kappa(B)$ be the least common multiple of $2$ and the positive
integers $n\geq 3$ such that

\begin{enumerate}

\item [(i)] $\Q (\zeta_n +\zeta_n^{-1})\subseteq F$, and

\item[(ii)] For all $\ell \mid \Norm_{F/\Q }(\mathfrak D)$, no prime ideal $\wp \mid \frac{\mathfrak D}
{(\mathfrak D, \mathfrak a_{\ell })}$ splits in $F(\zeta_n)$.
\end{enumerate}

\end{defn}

For instance, when $F=\Q $ or $F$ quadratic, $F\ne \Q(\sqrt{2})$, $\Q(\sqrt{5})$, we have $\kappa (B)\mid 12$. If $F=\Q(\sqrt{2})$, $\kappa (B)\mid 24$. If $F=\Q(\sqrt{5})$, $\kappa (B)\mid 60$.

For arbitrary $F$, we have $\kappa (B)\mid \text{lcm}_{\varphi(n)\mid g} (n)$, since $[\Q(\zeta_n):\Q ]=\varphi (n)=\sharp (\Z/n\Z )^*$. Also,
$8\nmid \kappa (B)$ provided $\sqrt{2}\not\in F$, because $\Q (\zeta_8+\zeta_8^{-1}) = \Q(\sqrt{2})$.

\begin{lema}\label{lcm}
$\kappa (A, \mathfrak N)\mid 2 \cdot \kappa (B)$.
\end{lema}

{\em Proof. } Let $\ell $ be a prime and $v$ a place of $K$ over
$\ell $. Since $[I_{\ell }: I_v]\leq 2$, it suffices to show that $\alpha^{\kappa (B)
}_{\mathfrak N}(I_{v})=\{ 1 \}$ for all places $v\nmid N$ of $K$.

Let $L_{\ell }/\Q_{\ell }$ be as in Lemma\,\ref{good} the extension constructed by Serre and Tate over which $A$
acquires good reduction. Let $\tilde
A/\F_{\ell }$ denote the closed fibre of the N\'eron model $\mathcal A$ of the
abelian variety $A\otimes L_{\ell }$.

Let $\bar{\F}_{\ell}$ be a fixed algebraic closure of $\F_{\ell }$.
Let $\tilde B=\End_{\bar{\F}_{\ell}}(\tilde A)\otimes \Q $ be the
algebra of endomorphisms of $\tilde A\otimes \bar{\F}_{\ell}$, in
which $B$ is embedded.

The group $I_v$ acts on $\mathcal A\times K_v$ via the Galois action on $L_{\ell }\cdot K_v$ over $K_v$. Since $B=\End^0_K(A)$, the action of $I_v$ commutes with that of $B$.

In addition, $I_v$ acts trivially on the residue field of $L_{\ell }\cdot K_v$ and thus specializes to an action on the closed fibre of $\mathcal A$ by algebraic automorphisms (cf.\,\cite[p.\,497]{SeTa}). By considering the action both of $I_v$ and $\tilde B$ on $V_{\mathfrak N}(A\times L_{\ell }\cdot K_v)=V_{\mathfrak N}(\tilde A)$ we obtain that $r_{\mathfrak N}(I_v)$ is a finite subgroup of the group
$G$ of units of an order in $C_{\tilde B}(B)$, the centralizer of
$B$ in $\tilde B$.

As in the proof of Proposition 5.2 in \cite{Mip}, $C_{\tilde B}(B)$ is either
$(i)$ a totally imaginary quadratic extension of $F$ which splits
$B$ or $(ii)$ $B(\ell )$.

In $(i)$, $A$ is isogenous over $\bar{\F}_{\ell}$ to
$A_0^2$, where $A_0/\bar{\F}_{\ell}$ is an abelian variety of
dimension $g/2$ with $\End_{\bar{\F}_{\ell}}(A_0)\supseteq C_{\tilde B}(B)$. In
$(ii)$ $A/\bar{\F}_{\ell}$ is supersingular: it is
isogenous to $g$ copies of a supersingular elliptic curve.

If $(i)$ holds, $G$ is the cyclic
group of $n$-th roots of unity in the CM-field $C_{\tilde B}(B)$ for some $n\geq 1$. As soon as $n>2$,
$C_{\tilde B}(B)=F(\zeta_n)$. Since it splits $B$, no prime $\wp \mid \mathfrak D$ splits in $F(\zeta_n)$.

Assume now that $(ii)$ holds. If $G\supseteq \langle \zeta_n\rangle $ for some $n>2$, then $F(\zeta_n)$ is a quadratic extension of $F$ which embeds in $B(\ell )$. We again deduce that no $\wp \mid \disc(B(\ell ))$ splits in $F(\zeta_n)$. If $\ell \nmid \Norm_{F/\Q }(\mathfrak D)$, this implies that
no $\wp \mid \mathfrak D\cdot \mathfrak a_{\ell }$ splits in $F(\zeta_n)$.
If $\ell \mid \Norm_{F/\Q }(\mathfrak D)$, we obtain that
no $\wp \mid \frac{\mathfrak D}{(\mathfrak D, \mathfrak a_{\ell })}$ splits in $F(\zeta_n)$.

By Lemma \ref{Borel}, there exists a surjective homomorphism from
$r_{\mathfrak N}(I_v)$ onto $\alpha_{\mathfrak N}(I_v)$. Since
$\alpha_{\mathfrak N}(I_v)\subseteq k^*$, this is a cyclic group.
The discussion above yields that $\alpha_{\mathfrak N}(I_v)^n=\{
1\}$ for some $n$ as in the statement. $\Box $

\vspace{0.2cm}

The above proof shows a slightly stronger fact: If $\ell \nmid \Norm_{F/\Q }(\mathfrak D)$ is a {\em supersingular prime} for $A$, then no $\wp \mid \mathfrak D\cdot \mathfrak a_{\ell }$ splits in $F(\zeta_n)$. We did not include this in the statement of Lemma \ref{lcm} because we shall apply it to {\em any} modular abelian variety $A/\Q $ with multiplication by $B$. One can not expect to find a prime $\ell \nmid \Norm_{F/\Q }(\mathfrak D)$ which is supersingular for all such $A$ simultaneously.

\begin{cor}\label{coco}
Assume $K=\Q(\sqrt{-N})$ and $4\mid \kappa (B)$. Then
$\kappa (A, \mathfrak N) \mid \kappa (B)$.
\end{cor}

{\em Proof.} In the above proof we showed that $\alpha^{\kappa (B)
}_{\mathfrak N}(I_{v})=\{ 1 \}$ for all places $v\nmid N$ of $K$. Since $K=\Q(\sqrt{-N})$, we have $I_v = I_{\ell }$ for all these places of $K$. As $4\mid \kappa (B)$, it follows from the definition of $\kappa (A, \mathfrak N)$ that it divides $\kappa (B)$.  $\Box $

\begin{prop}\label{caracter}  For any prime $\ell \ne N$, $\alpha_{\mathfrak N}(\varphi_{\ell }^{\kappa (A, \mathfrak N)}) =  {\ell }^{\kappa (A, \mathfrak N)/2}$.
\end{prop}

{\em Proof. } Let $\kappa = \kappa (A, \mathfrak N)$. Let us first show that $\alpha_{\mathfrak
N}(\varphi_{\ell }^{\kappa}) = \varrho ^{\kappa }_N(\frac{1}{\ell })
\in \F_{N}^*$. The image by $\om $ of the tuple $$(\frac{1}{\ell },
..., \frac{1}{\ell }, \stackrel{\stackrel{{\ell }}{\smallsmile}}{1},
\frac{1}{\ell }, ..., \frac{1}{\ell }) \in \prod_{p}\Z_{p}^*$$ is an
element $\tilde\varphi_{\ell }\in \GQ^{ab}$ which reduces to the
Frobenius automorphism $\Fr_{\ell} \in \Gal(\bar \F_{\ell }/\F_{\ell
})$. Hence $\varphi_{\ell }\circ \tilde \varphi_{\ell }^{-1}\in
I_{\ell }$ and $\alpha_{\mathfrak N}(\varphi_{\ell }^{\kappa }) =
\alpha_{\mathfrak N}(\tilde\varphi_{\ell}^{\kappa })$ because
$\alpha_{\mathfrak N}^{\kappa }(I_{\ell }) = \{ 1\}$. It follows
that $\alpha^{\kappa}_{\mathfrak N} (\varphi_{\ell }) = \varrho
^{\kappa }_{\ell }(1)\cdot \varrho^{\kappa }_N (\frac{1}{\ell
})\cdot \prod_{p\ne N, \ell }\varrho ^{\kappa }_{p}(\frac{1}{\ell })
= \varrho^{\kappa }_N (\frac{1}{\ell })$.

Let $\bar \chi_N : \GQ \, \lra \, \F^*_N$ be the reduction of the
cyclotomic character mod $N$. By \cite[Prop. $3$, $8$]{Se} we
have $\bar \chi_{N |I_N^{ab}}\cdot  \om_N (x)= \frac{1}x$ mod $N$
for all $x\in \Z_N^*$.

Since $\F_N^*$ has order $N-1$, we have $$\varrho _N (x)= x^{-c }
\mbox{ mod } N \mbox{ for some }0\le c<N-1.$$ By Lemma \ref{Borel},
since $\det (r_{\mathfrak N}) =\chi_N$ we have $\frac{1}x = $ $\bar \chi_N(\om_N(x))
= $ $\det \bar r_{\mathfrak N}(\om_N(x)) = $ $\psi(\om_N(x))\cdot
N_{k/k_0} (\alpha_{\mathfrak N}(\om_N(x)) = $ $\psi(\om_N(x))
N_{k/k_0}(\varrho_N(x)).$

Since $\varrho_N(x)\in \F_N^*$, $N_{k/k_0}(\varrho_N(x)) =
\varrho^2_N(x)$ and it follows that $\varrho^2_N(x) = x^{-2 c} = \pm
x^{-1}$. Hence $2 c \equiv 1 \,\mathrm{ or }\,\frac{N+1}2$
mod $N-1$. Since $N-1$ is even, we deduce that $N\equiv 3$ mod $4$
and $c = \frac{N+1}4$ or $\frac{3N-1}4$ mod $N-1$. Thus
$\varrho^{\kappa }_N(x) \equiv x^{-\frac{(N+1) \kappa }4}$ or
$x^{-\frac{(3N-1) \kappa }4}$ mod $N$. Since $4\mid \kappa$, $\frac{(N+1) \kappa }4 \equiv
\frac{(3N-1) \kappa }4 \equiv \kappa/2 $ mod $N-1$ and
$\alpha_{\mathfrak N}(\varphi_{\ell }^{\kappa }) = \ell ^{\kappa/2} \in \F_{N}^*$. Consequently,
$\alpha_{\mathfrak N}(\varphi_{\ell }^{\kappa}) = {\ell }^{\kappa/2}$. $\Box $

\vspace{0.3cm}
During the proof of the above proposition we proved in passing the following.

\begin{cor}\label{N}
$N\equiv 3$ mod $4$.
\end{cor}

\begin{cor}\label{kappa} If $\mathrm{ord}\,_2(\kappa (A, \mathfrak N) )\leq \mathrm{ord}\,_2(\nu +1)$ then $K=\Q(\sqrt{-N})$.
\end{cor}

{\em Proof. } Let $\kappa = \kappa (A, \mathfrak N)$. For each $\ell \ne N$ we have $\ell \equiv \psi
(\Fr_{\ell }) \alpha_{\mathfrak N} (\Fr_{\ell })^{\nu+1}$ mod
$\mathfrak N$. Set $t=\kappa /2^{\mathrm{ord}_2\,(\kappa )}$. Then
$\kappa \mid t (\nu +1)$ and $\ell ^t\equiv \psi^t(\Fr_{\ell })
\alpha_{\mathfrak N}^{t (\nu+1)}(\Fr_{\ell })$. Since $t$ is odd and
$\psi $ is a quadratic character, it follows from Proposition
\ref{caracter} that $\psi(\Fr_{\ell }) \equiv \ell^{t
(\nu-1)/2}$ mod $\mathfrak N$. As $\psi(\Fr_{\ell }), \ell \in \Z $ and
$\mathfrak N=\mathfrak N_0\cdot R$, this is equivalent to the congruence
$\psi(\Fr_{\ell }) \equiv \ell^{t(\nu-1)/2}$ mod $\mathfrak N_0$.

Since the residue field of $\mathfrak
N_0$ is $k_0$, we obtain that
$\psi(\Fr_{\ell }) = (\frac{\ell }{k_0})$ for all $\ell \ne
N$. If $f$ is even, then $(\frac{\ell }{k_0})=1$ for all $\ell $ and
$\psi $ is trivial, a contradiction. Hence $f$ is odd and
$(\frac{\ell }{k_0}) =(\frac{\ell }{N}) = (\frac{-N}{\ell })$. It follows that
$\psi = (\frac{-N}{\cdot })$ and
$K=\Q(\sqrt{-N})$. $\Box $

\vspace{0.2cm}

One application of Corollary \ref{kappa} is in Theorem \ref{Sh}
$(ii)$. It can also be used to show in many instances that an
analogue of (\ref{SY}) in the Introduction holds for higher degree
fields $F$.

Let for example $(\cO_{\mathfrak D}, R_{F(\sqrt{m})}, K)$ be a
modular triplet such that $\sqrt{2}\not \in F$ and $m\mid \mathfrak
D=\mathfrak N_1\cdot ... \cdot \mathfrak N_{2 r}$, above integer
primes $N_i$. Order them in some way that $\mathfrak N_i\mid 2 m$
for $1\leq i \leq s$ and $\mathfrak N_i\nmid 2 m$ for $i>s$. Assume
$\mathfrak N_{2 r-1}$ and $\mathfrak N_{2 r}$ can be chosen such
that $N_{2 r-1}, N_{2 r}$ are odd, $N_{2 r-1}\ne N_{2 r}$ and
$f(\mathfrak N_{2 r}/N_{2 r})$ is odd.

If either $(a)$ $N_i, N_j\equiv 1$ mod $4$ for some $1\leq i\ne
j\leq s$ or $(b)$ $N_{2 r}\equiv 7$ mod $8$, then $s\geq 2 r-1$.

Indeed, if $s\leq 2 r-2$, it would follow from (\ref{ad}) in the
Introduction that $N_{2r -1} N_{2 r}\mid \disc (K)$. By Definition
\ref{B}, $4\nmid \kappa (B)$ in $(a)$, $8\nmid \kappa (B)$ in $(b)$.
By Lemma \ref{lcm}, $8\nmid \kappa (A, \mathfrak N_{2 r})$ in $(a)$,
$16\nmid \kappa (A, \mathfrak N_{2 r})$ in $(b)$. Corollary
\ref{kappa} applied to $\mathfrak N_{2 r}$ would imply $K=\Q
(\sqrt{-N_{2 r}})$, a contradiction.

\begin{lema}\label{fites}  The characteristic polynomial of $\varphi_{\ell }$ satisfies
$P_{\ell }(T) \,\equiv \, T^2 - (\alpha_{\mathfrak N} (\varphi_{\ell
})+ {\ell } \alpha_{\mathfrak N} (\varphi^{-1}_{\ell }))\cdot T +
{\ell } \,\in k[T]$.

\end{lema}

{\em Proof. } It follows from Lemma \ref{Borel} that for any $\sigma
\in \GQ $, the characteristic polynomial $P_{\sigma }(T)\in R[T]$
satisfies
$$
P_{\sigma } \,\mathrm{mod}\, \mathfrak N \equiv
T^2 - (\alpha_{\mathfrak N} (\sigma )+\psi (\varphi_{\ell }) \alpha_{\mathfrak N} (\sigma
)^{\nu })\cdot T + \psi (\varphi_{\ell }) \mathrm{Nm}(\alpha_{\mathfrak N} (\sigma )).
$$

Besides, as we mentioned at the end of Section $2$, we have $P_{\ell }=P_{\varphi_{\ell }} = T^2-a_{\ell } T + \ell $.
Hence $\psi (\varphi_{\ell }) \alpha_{\mathfrak N} (\varphi_{\ell })^{\nu } \alpha_{\mathfrak N}
(\varphi_{\ell }) = {\ell } \in k$ and we deduce that $a_{\ell
}\equiv \alpha_{\mathfrak N} (\varphi_{\ell }) + \psi (\varphi_{\ell })
\alpha_{\mathfrak N} (\varphi_{\ell })^{\nu }\equiv \alpha_{\mathfrak N} (\varphi_{\ell })+ {\ell
} \alpha_{\mathfrak N} (\varphi^{-1}_{\ell })$ mod $\mathfrak N$. $\Box $

\section{Proof of the main results}\label{last}

Keep the notations of the previous section. Let $\lambda_{\pm} =
\sqrt{2 \ell \pm \sqrt{3}}$.

\begin{teor}\label{MaIn} Let $(\cO, R, K)$ be a modular triplet. Assume
$\zeta_n+\zeta_n^{-1}\not\in F$ for any $n$-th primitive root of
$1$, $n\ne 1, 2, 3, 4, 6$, and that $\cO $ is locally maximal at
some ideal $\mathfrak N_0\mid \disc(B)$, $\mathfrak N_0\nmid 2 m$.
Then, for any prime $\ell $ which splits or ramifies in $K$, either
$\mathfrak N_0\in \mathcal N_{\ell }$ or

\begin{itemize}

\item $(\frac{-\ell }{\wp }) \ne 1$ for all $\wp \mid \mathfrak D$, or

\item $\sqrt{2 \ell }\in F$ and $(\frac{-1 }{\wp }) \ne 1$ for all $\wp \mid \mathfrak D$, $\wp \nmid \ell $, or

\item $\sqrt{\ell }\text{ or }\sqrt{3 \ell }\in F$ and $(\frac{-3 }{\wp }) \ne 1$ for all $\wp \mid \mathfrak D$, $\wp \nmid \ell $, or

\item $\lambda_{\pm } \sqrt{\ell }\in F$ and $(\frac{-1\mp 4\sqrt{3}/7}{\wp }) \ne 1$ for all $\wp \mid \mathfrak D$, $\wp \nmid \ell $.

\end{itemize}

\end{teor}

{\em Proof. } Let $A$ be a modular abelian variety over $\Q $ such that $\End_{\Q }(A)=R$ and $\End_K(A)\simeq \cO $.

By our assumption on $\mathfrak N_0$, it remains inert in $R$ and we
let $\mathfrak N = \mathfrak N_0 R$ (cf.\,Section \ref{Ass}). Let
$\bar{\mathfrak N}$ be an ideal of $\qbar $ above $\mathfrak N_0$.
Choose a prime $\ell $ such that $\mathfrak N_0\not \in \mathcal
N_{\ell }$ and $(\frac{K}{\ell })\ne -1$. Then $\mathfrak N_0\nmid
\ell $, $\varphi_{\ell }\in G_K$ and $a_{\ell }\in R_F$.

By Lemma \ref{fites}, $\alpha_{\mathfrak N} (\varphi_{\ell })+ {\ell
} \alpha_{\mathfrak N} (\varphi^{-1}_{\ell }) = \sqrt{{\ell }}\cdot
(\frac{\alpha_{\mathfrak N} (\varphi_{\ell })}{\sqrt{{\ell }}} +
(\frac{\alpha_{\mathfrak N} (\varphi_{\ell })}{\sqrt{{\ell
}}})^{-1}) \equiv a_{\ell }$ mod $\bar{\mathfrak N}$. Here,
$\zeta=\frac{\alpha_{\mathfrak N} (\varphi_{\ell })}{\sqrt{{\ell
}}}$ is a $\kappa (A, \mathfrak N)$-th root of $1$ by Proposition
\ref{caracter}.

Since $\zeta_n+\zeta_n^{-1}\not\in F$ for any $n\ne 1, 2, 3, 4, 6$,
we have $\kappa (B)\mid 12$. Hence $\kappa (A, \mathfrak N)\mid 24$
by Lemma \ref{lcm}.

Computing the possible values of $\sqrt{\ell }(\zeta +\zeta^{-1})$
for all $24$-th roots of $1$, we obtain that $a_{\ell } \equiv 0,
\sqrt{\ell },\, \sqrt{ 2 \ell}, \sqrt{3 {\ell }},\,2\sqrt{ \ell }$
or $\lambda_{\pm}\cdot \sqrt{\ell }$ mod $\bar {\mathfrak N}$. In
other words, $\mathfrak N_0\mid a_{\ell }^2 - s \ell $ for some $s=
0, 1, 2, 3, 4$ or $\mathfrak N_0\mid a_{\ell }^4  - 4 a_{\ell }^2
\ell + \ell^2$.

As we already mentioned, $|\,\tau ({a_{\ell }})|\le 2 \sqrt{{\ell
}}$ for any $\tau : F\hookrightarrow \R $. It follows from the very
definition of $\mathcal N_{\ell }$ that the above congruence is in
fact an identity: $a_{\ell }= 0, \sqrt{\ell },\, \sqrt{ 2 \ell},
\sqrt{3 {\ell }},\,2\sqrt{ \ell }$ or $\lambda_{\pm}\cdot \sqrt{\ell
}$.

Since $(\frac{K}{\ell })\ne -1$, the residue field of $K$ at $\ell $
is $\F_{\ell }$. By reducing the endomorphisms of $A$ mod $\ell $ it
follows that $B=\End_K(A)\otimes \Q $ embeds in $\End_{\F_{\ell
}}(\tilde A)\otimes \Q $.

According to \cite[Theorem 2 (a)]{Ta}, \cite[Theorem 1 (1)]{Ta2},
$\End_{\F_{\ell }}(\tilde A)\otimes \Q $ is isomorphic to $\M_{2
r}(\tilde F)$, where $\tilde F$ is the splitting field of
$\prod_{\tau: F\hookrightarrow \R }(T^2-a^{\tau }_{\ell } T+\ell )$.
Here, $r\mid g=\dim(A)$ and $[\tilde F:\Q ]=g/r$.

From the above values for $a_{\ell }$  we have (in the same order)
$\tilde F = \Q(\sqrt{-\ell }),$ $\Q(\sqrt{\ell},$ $\sqrt{-3}),$
$\Q(\sqrt{2 \ell},$ $\sqrt{-1}),$ $\Q(\sqrt{3 \ell}, \sqrt{-3})$,
$\Q (\sqrt{\ell })$ or $\Q(\lambda_{\pm} \cdot \sqrt{\ell },
\sqrt{-1\mp 4 \sqrt{3}/7})$.

Case $\tilde F=\Q(\sqrt{\ell })$ can not arise: we would have
$\End_{\F_{\ell }}(A)\otimes \Q \simeq \M_g(\Q(\sqrt{\ell }))$ and
this would say that $A$ is isogenous to $g$ copies of an elliptic
curve $E/\F_{\ell }$ with $\End_{\F_{\ell }}(E)\otimes \Q \simeq
\Q(\sqrt{\ell })$, which is not possible.

Let us show that $F\cdot \tilde F$ splits $B$ locally at all places
$\wp \nmid \ell $, or equivalently, no prime $\wp \mid \mathfrak D$,
$\wp \nmid \ell$, splits in the quadratic extension $F\cdot \tilde
F$ of $F$.

The Tate module $V_p(\tilde A)$ is a $F\cdot \tilde F\otimes
\Q_p$-vector space of rank $2$. Since $B\subset \End_{\F_{\ell
}}(\tilde A)\otimes \Q $ and $\tilde F$ is the center of
$\End_{\F_{\ell }}(\tilde A)\otimes \Q $, we conclude that $B$ acts
$F\cdot \tilde F\otimes \Q_p$-linearly on  $V_p(\tilde A)$.

Thus $B\subset \M_2(F\cdot \tilde F\otimes \Q_p)$, which says that
$F\cdot \tilde F$ splits $B$ locally at $\wp $.

Finally, when $\tilde F = \Q (\sqrt{-\ell })$, it is also clear that
$F\cdot \tilde F = F(\sqrt{-\ell })$ also splits $B$ at the places
$\wp \mid \ell $, because these prime ideals ramify in $F\cdot
\tilde F$. This finishes the proof of the theorem. $\Box $

\vspace{0.2cm}

Theorem \ref{main} $(iii)$ now follows automatically by applying
Theorem \ref{MaIn} to each prime ideal $\mathfrak N \mid \mathfrak
D$, $\mathfrak N\nmid 2 m$. Note that, in this case, since we assume
that $\sqrt{\ell }$, $\sqrt{2 \ell }$, $\sqrt{3 \ell }$,
$\sqrt{2\ell \pm \sqrt{3} \ell }\,\not \in F$, it follows that
$a_{\ell } = 0$ and $\tilde F = \Q(\sqrt{-\ell })$.

\begin{remark} As we saw in Proposition \ref{unr},
Theorem \ref{Sh} and Corollary \ref{kappa}, we can often claim that
$K=\Q (\sqrt{-N})$. When this is the case, last item in the
statement of Theorem \ref{MaIn} can be removed and one can replace
the sets $\mathcal N_{\ell }$ by the smaller sets $\mathcal
N^0_{\ell } = \{ \mathfrak N_0 \mid \ell \text{ or } a^2- s \ell
\}$, with $a$ and $s$ as in Definition \ref{Nl}. This is because in
the above proof we have $\kappa (A, \mathfrak N)\mid 12$, thanks to
Corollary \ref{coco}.
\end{remark}

\vspace{0.2cm}

{\em Proof of Theorem \ref{Sh}. } Since the hypothesis of Theorem
\ref{main} $(ii)$ are satisfied, we obtain that $N\equiv 3$ mod $4$.

Assume first that $M\equiv 3$ mod $4$. Proposition \ref{unr} applies
and asserts that $K\simeq \Q(\sqrt{-d})$, where $d=M$ or $N$. Since
$B\simeq (\frac{-d, M }{\Q })$ we deduce that $d=N$ and
$(\frac{-N}M)=-1$. Finally, for any odd prime $\ell $ such that
$(\frac{\ell }{N}) = 1$ and $N\not \in \mathcal N_{\ell }$, Theorem
\ref{MaIn} together with Gauss' reciprocity law asserts that
$(\frac{-\ell }M) = -1$.

Assume now that $M\equiv 1$ mod $4$. Proposition \ref{unr} now
implies that $K\simeq \Q(\sqrt{-d})$ for $d=M N$ or $N$.  In any
case, it again follows as above that $(\frac{-N}M)=-1$.

If $d=N$, we similarly have $(\frac{-\ell }M) = -1$ for any prime
$\ell $ such that $(\frac{\ell }{N}) = 1$ and $N\not \in \mathcal
N_{\ell }$, even for $\ell =2$.

If $d=M N$, $N\equiv 3$ mod $8$ by Lemma \ref{lcm} and Corollary
\ref{kappa}. Let $\ell \ne 2$ be a prime such that $N\not \in
N_{\ell }$ and $(\frac{-M N}{\ell }) = 1$. Theorem \ref{MaIn}
produces the congruences $(\frac{-\ell }N) = (\frac{-\ell }M) = -1$,
which Gauss' reciprocity law shows to be incompatible with
$(\frac{-M N}{\ell }) = 1$. $\Box $

\end{document}